\documentclass{amsart}
\usepackage{amsmath,amscd,amsthm,amsfonts,amssymb}

\theoremstyle{plain}

\theoremstyle{definition}

\numberwithin{equation}{section}

\def \theo-intro#1#2 {\vskip .25cm\noindent{\bf Theorem #1\ }{\it #2}}

\def\eps{{\varepsilon}}

\def\X{\Xi}
\def\x{\xi}

\def\lpar{\left (}
\def\rpar{\right )}
\def\lmpar{\left \{ }
\def\rmpar{\right \}}
\def\lkpar{\left [}
\def\rkpar{\right ]}

\def\labs{\left |}
\def\rabs{\right |}
\def\babs#1{\labs {#1} \rabs}

\def\dx{\, dx}

\def\Cbf{{\mathbf C}}

\def\Ccal{{\mathcal C}}
\def\Fcal{{\mathcal F}}
\def\Gcal{{\mathcal G}}

\def\Qcal{{\mathcal Q}}
\def\Rbf{{\mathbf R}}
\def\Rcal{{\mathcal R}}
\def\Rn{{\mathbf R \sp n}}
\def\Bn{{B \sp n}}
\def\Nbf{{\mathbf N}}

\def\Tbf{{\mathbf T}}
\def\Tn{{\mathbf T \sp n}}
\def\Zn{{\mathbf Z \sp n}}
\def\Wcal{{\mathcal W}}
\def\Xcal{{\mathcal X}}
\def\Ycal{{\mathcal Y}}
\def\Zbf{{\mathbf Z}}
\def\Zcal{{\mathcal Z}}
\def\Ycal{{\mathcal Y}}

\def\Co#1#2{#1 {\Ccal \sb 0} \lpar \Rbf \sp #2 \rpar}
\def\co#1#2{#1 {c \sb 0} \lpar {#2} \rpar}

\def\Lsp#1#2{{L \sp {#2}\lpar\Rbf \sp {#1}\rpar}}

\def\lsp#1#2{{\ell \sp {#2} \lpar {#1} \rpar}}

\def\nm#1#2#3{\left\|#1\right\| \sp {#2} \sb {#3}}

\def\ph{\varphi}

\def\xx{x \xi}

\def\axi{|\x|}

\def\la{\lambda}
\def\supp{\operatorname{supp}}

\begin{document}

\baselineskip 22pt \larger

\allowdisplaybreaks

\title
{On the Finite Dimensionality of Spaces of \\
 Absolutely Convergent Fourier Transforms} \author
{Bj\"orn G.\ Walther} \date
{\today} \keywords
{Absolutely convergent Fourier transforms, reflexivity, weakly sequential completeness} \subjclass
[2000]{42A38, 42A55, 42B10, 42B35, 46E15} \address
{Stockholm University, Department of Mathematics, SE -- 106 91 STOCKHOLM, Sweden.} \email
{walther@math.su.se}

\begin{abstract}
 We extend the result of K.\ Karlander \cite {Karlander_1997} regarding finite dimensionality of spaces of absolutely convergent Fourier transforms. \end{abstract} \maketitle \section
{Introduction} \subsection {} \label {sbsc_1.1}
The space of absolutely convergent Fourier transforms is dense in the space of continuous functions which vanish at infinity. See Segal \cite {Segal_1950}. Hence we may claim that sufficiently many functions are absolutely convergent Fourier transforms. On the other hand, the space of absolutely convergent Fourier transforms either coincides with the space of continuous functions which vanish at infinity or is of the first category in that space. This follows from a theorem on bounded linear mappings (cf.\ Banach \cite [th\'eor\`eme 3 p.\ 38] {Banach_1932}, \cite [theorem 3 p.\ 24] {Banach_1987}) applied to the {\it Fourier transformation}. The image of that bounded linear mapping being of the first category manifests that few functions are absolutely convergent Fourier tranforms. \subsection {}
A further indication of the absolutely convergent Fourier transforms being few is the following result: \subsubsection {\bf
 Theorem {\rm (K. Karlander \cite {Karlander_1997})}} \label {thm_1.2.1} {\it
 Let $\Ycal$ be a closed subspace of the space of continuous functions on the real line which vanish at infinity. \linebreak Assume that all elements in $\Ycal$ are absolutely convergent Fourier transforms. If in addition $\Ycal$ is reflexive, then $\Ycal$ is of finite dimension.} \subsection {
 The Purpose of this Paper} Our purpose here is to discuss on and extend theorem \ref {thm_1.2.1}. See theorem \ref {thm_4.1} on page \pageref {thm_4.1}. In the course of proving this theorem we extract two lemmata which we beleive are of independent interest. See \S\S \ref {lma_6.1} and \ref {lma_6.5} on page \pageref {lma_6.1} and \pageref {lma_6.5} respectively. In extending the given result we replace $\Lsp {} 1$ by $\Lsp n 1$ and $L \sp 1 \lpar \Tn \rpar$ respectively. The proofs follow the ideas in \cite {Karlander_1997} closely. \subsection {
 The Plan of this Paper}
 In \S \ref {sec_2} we introduce the notation needed for this paper. In \S \ref {sec_3} starting on page \pageref {sec_3} we discuss on some earlier and related results. In \S \ref {sec_4} on page \pageref {sec_4} we state our main result. In \S \ref {sec_5} starting on page \pageref {sec_5} we cite results needed for our proofs, and in \S \ref {sec_6} starting on page \pageref {sec_6} we provide specific tools in the form of some lemmata. A proof of our main result is given in \S \ref {sec_7} starting on page \pageref {sec_7}, and in \S \ref {sec_8} starting on page \pageref {sec_8} we provide some examples of closed non-reflexive subspaces of infinite dimension. \section {
 Notation} \label {sec_2} \subsection {
 Banach Spaces} By $\Wcal$ we denote a Banach space with dual space $\Wcal \sp *$. With few exceptions all function spaces are Banach spaces. \subsubsection {
 Lebesgue Spaces} We define $L \sp 1 (S, \Sigma, \mu)$ in the usual way. In most cases $S$ will be either $\Rn$ or $\Tn$ endowed with Lebesgue measure. We will then write $\Lsp n 1$ and $L \sp 1 \lpar \Tn \rpar$ respectively. \subsubsection {
 Spaces of Convergent Sequences and of Summable Sequences} By $\co {} M$ we denote the set of all {\it functions} $a$ on $M$ such that $a(m) = 0$ for all but finitely many choices of $m$. By $\co \overline M$ we denote the completion of $\co {} M$ under the norm given by
 \begin{equation*}
 \nm a {} {\co \overline M} \, = \, \sup \sb {m \in M} \babs {a(m)}. \end{equation*}
 We define $\lsp M 1$ in the usual way. \subsubsection {
 Spaces of Continuous Functions} By $\Tbf$ we denote the half-open \linebreak interval $[-\pi,\pi[$. The space $\Ccal \lpar \Tbf \rpar$ is the space of functions $g$ continuous on $\Tbf$ such that $g(\pi-)$ exists and equals $g(-\pi)$. This requirement allows us to regard $g$ as a continuous function on the unit circle in $\Cbf$. Thus
 \begin{equation*}
 \nm g {} {\Ccal \lpar \Tbf \rpar} \, = \, \sup \sb {x \in \Tbf} \babs {g(x)} \end{equation*}
 exists for every $g \in \Ccal \lpar \Tbf \rpar$ and this quantity defines a norm on $\Ccal \lpar \Tbf \rpar$. \par
 Spaces of continuous functions on other compact sets will occur and their definitions and norms are obvious. \par
 The space $\Co {} n$ is the space of functions $F$ continuous on $\Rn$ such that $\supp F$ is compact. By $\Co \overline n$ we denote the completion of $\Co {} n$ under the norm given by
 \begin{equation*}
 \nm F {} {\Co \overline n} \, = \, \sup \sb {\x \in \Rn} \babs {F(\x)}. \end{equation*}
 $X$ will denote a locally compact abelian group with dual group $\X$. To obtain $\overline {\Ccal \sb 0}(\X)$ we apply the previous definition regarding completion. Observe that we have $\overline {\Ccal \sb 0} \lpar \Zbf \sp n \rpar = \co \overline {\Zbf \sp n}$. \subsection {
 The Fourier Transformation} If the action of $\x \in \X$ on $x \in X$ is written $e \sp {i\xx}$ and if the Haar measure on $X$ is written
 \begin {equation*}
 E \longmapsto \int \sb E \dx \end{equation*}
 then for $f \in L \sp 1 (X)$ the integral
 \begin{equation*}
 \lkpar \Fcal f \rkpar \lpar \x \rpar \, = \, \int \sb X \overline {e \sp {i\xx}} f(x) \dx \end{equation*}
 is absolutely convergent for all $\x \in \X$ and the function $\Fcal f$, the {\it Fourier transform} of $f$, is in $\overline {\Ccal \sb 0}(\X)$. The mapping $\Fcal$, the {\it Fourier transformation,} is bounded, linear and injective $L \sp 1 (X) \longrightarrow \overline {\Ccal \sb 0} (\X)$. \par
 By $\lkpar \Fcal L \sp 1 \rkpar (\X)$ we denote the {\it range} of the Fourier transformation $\Fcal$. We endow this space with the topology induced by $\overline {\Ccal \sb 0} (\X)$. \subsection {
 Auxiliary Notation} By $\overline {\Bn}$ we denote the closed unit ball of $\Rn$. Numbers denoted by $C$ may be different at each occurrence. At several occasions the notation $C \sb {p.q}$ will be used, where $p$ and $q$ are positive integers. In this notation the double subscript $p.q$ refers to inequality $(p.q)$ from which the existence of $C \sb {p.q}$ follows. \section {
 Earlier Results} \label {sec_3} \subsection {
 Density} As already mentioned in \S \ref {sbsc_1.1} on page \pageref {sbsc_1.1}, the range $[\Fcal L \sp 1](\X)$ is dense in the codomain $\overline {\Ccal \sb 0}(\X)$ of $\Fcal$. See Segal \cite [lemma 2 p.\ 158] {Segal_1950}. This holds for any $X$. The density may be proved indirectly using the separation \linebreak theorem of Hahn and Banach, the representation theorem of Riesz and Markov and the theorem of Fubini on the change of order of integration. The density result extends to the transformation of Fourier and Stieltjes. See Hewitt \cite [theorem 1.3 p.\ 664] {Hewitt_1953}. \subsection {
 Category and Cardinality} The space $[\Fcal L \sp 1](\X)$ is of the first \linebreak category in the codomain unless it coincides with the codomain. This is a general result for bounded linear mappings. Moreover, if $\lkpar \Fcal L \sp 1 \rkpar (\X)$ coincides \linebreak with the codomain, then $X$ is finite. See Segal \cite [p.\ 157] {Segal_1950}. Cf.\ also Edwards \cite {Edwards_1954}, Rajagopalan \cite {Rajagopalan_1964}, Friedberg \cite {Friedberg_1971}, Graham \cite {Graham_1973} and Basit \cite {Basit_1981}. A similar \linebreak result is valid for the transformation of Fourier and Stieltjes. See Hewitt \cite [theorem 1.2 p.\ 664] {Hewitt_1953}. Cf.\ also Edwards \cite {Edwards_1954}. \subsection {
 Finite Dimensionality} The result in \cite {Karlander_1997} extended in theorem \ref {thm_4.1} on page \pageref {thm_4.1} is of the following nature: On $\Ycal$, a closed subspace of the codomain of $\Fcal$, we impose additional conditions. The conslusion is that $\Ycal$ is of finite dimension. It is hence relevant to relate our results to earlier results which are of a similar pattern. \subsubsection {\bf
 Definition} \label {def_3.3.1} A sequence $(x \sb m) \sb {m \in \Zbf \sb +}$ is {\it weakly Cauchy} in $\Wcal$ if the sequence of complex numbers $(\langle x \sb m , x \sp * \rangle) \sb {m \in \Zbf \sb +}$ is Cauchy for every choice of $x \sp * \in \Wcal \sp *$. The space $\Wcal$ is {\it weakly sequentially complete} if every weakly Cauchy sequence in $\Wcal$ converges weakly. \subsubsection {\bf
 Theorem {\rm (Rajagopalan \cite [theorem 2 p.\ 87] {Rajagopalan_1964} and Sakai \cite [proposition 2 p.\ 661] {Sakai_1964}.)}} \label {thm_3.3.2} {\it
 Let $\Ycal$ be a weakly sequentially complete $C \sp *$-algebra. Then $\Ycal$ is of finite dimension.} \subsubsection {\bf
 Theorem {\rm (Cf.\ Albiac, Kalton \cite [corollary 2.3.8 p.\ 38] {Albiac_Kalton_2006}.)}} {\it
 Let $\Ycal$ be a Banach space such that every weakly convergent sequence is convergent. If $\Ycal$ is reflexive then it is of finite dimension.} \subsection {
 Necessary Conditions and Sufficient Conditions} In the survey paper by Liflyand, Samko and Trigub \cite {Liflyand_Samko_Trigub_2012} results on representation of a function as an absolutely convergent Fourier transform are collected. It contains an extensive list of references and also a table of functions which are absolutely convergent Fourier transforms. \section{
 Main Result} \label {sec_4} \subsection {
 Theorem} \label {thm_4.1} {\it
 Let $\Ycal$ be a closed subspace of $\overline {\Ccal \sb 0}(\X)$ such that $\Ycal$ is a subset of $\lkpar \Fcal L \sp 1 \rkpar (\X)$ where $X = \Rn$ or $X = \Tn$. If $\Ycal$ is reflexive then it is of finite dimension.} \par \ \par \noindent
 The theorem is proved in \S\S \ref {sbsc_7.1} and \ref {sbsc_7.2} on page \pageref {sbsc_7.1} and \pageref {sbsc_7.2} respectively. \section {
 Preparation} \label {sec_5} \subsection {}
 We refer to the definition of weakly sequential completeness of $\Wcal$ in \S \ref {def_3.3.1} on this page. \subsubsection {} \label {sbsbsc_5.1.1}
 The space $L \sp 1 (a,b)$ is weakly sequentially complete. This is a theorem of Steinhaus. Cf.\ \cite {Steinhaus_1919}. More generally, the space $L \sp 1 (S, \Sigma, \mu)$ for a positive measure space $(S, \Sigma, \mu)$ is weakly sequentially complete. Cf.\ e.g.\ Dunford, Schwartz \cite [theorem IV.8.6 p.\ 290] {Dunford_Schwartz_1958}. \subsubsection {} \label {sbsbsc_5.1.2}
 If $\lmpar e \sb m : m \in \Zbf \sb + \rmpar$ is the canonical basis in $\co \overline {\Zbf \sb +}$ then $(e \sb 1 + e \sb 2 + \dots + e \sb m) \sb {m \in \Zbf \sb +}$ is weakly Cauchy but does not converge weakly. Hence the space $\co \overline {\Zbf \sb +}$ is not weakly sequentially complete. Cf.\ e.g.\ Albiac, Kalton \cite [example 2.3.11 p.\ 38] {Albiac_Kalton_2006}. \subsection {
 Lemma} \label {lma_5.2} {\it
 Let $\Wcal \sb 1$ be a closed subspace of $\Wcal$. If $\Wcal$ is weakly sequentially complete, then $\Wcal \sb 1$ is weakly sequentially complete.} \par \ \par \noindent
 For the sake of completeness we provide a proof. \subsection* {\it
 Proof:} We need to show that every weakly Cauchy sequence in $\Wcal \sb 1$ converges weakly to an element $\Wcal \sb 1$. \par
 Assume that the sequence $(x \sb m) \sb {m \in {\Zbf \sb +}}$ in $\Wcal \sb 1$ is weakly Cauchy. Then the same sequence in $\Wcal$ is weakly Cauchy and hence it has a weak limit $x \in \Wcal$. If $x \notin \Wcal \sb 1$ then according to the separation theorem of Hahn and Banach there is $x \sp * \in \Wcal \sp *$ such that $\langle \Wcal \sb 1, x \sp * \rangle = \{ 0 \}$ and $\langle x, x \sp * \rangle = 1$. This contradicts $x$ being the weak limit given above. \subsection {
 Theorem \rm (Pe\l{}czy\'nski \cite [lemma 2 p.\ 214] {Pelczynski_1960}. Cf.\ also Lindenstrauss, Tzafriri \cite [proposition 2.a.2 p.\ 53] {Lindenstrauss_Tzafriri_1977}.)} {\it
 Let $\Ycal$ be a closed subspace of $\lsp {\Zbf \sb +} 1$ of infinite dimension. Then $\Ycal$ contains a subspace $\Zcal$ which is isomorphic to $\lsp {\Zbf \sb +} 1$.} \subsection {
 Lemma} {\it
 A closed subspace of a reflexive space is reflexive.} \subsection* {\it
 Comment on the proof:} A proof using the separation theorem of Hahn and Banach and the theorem of Banach and Alaoglu is outlined in Rudin \cite [exercise 1 p.\ 111] {Rudin_1991}. \subsection {
 Corollary} \label {cor_5.5} {\it
 Let $\Ycal$ be a closed subspace of $\lsp {\Zbf \sb +} 1 $ of infinite dimension. Then $\Ycal$ is not reflexive.} \subsection {
 Theorem \rm (Cf.\ Katznelson \cite [\S 1.4 p.\ 108] {Katznelson_1968}, \cite [\S 1.4 p.\ 137] {Katznelson_2004} and Zygmund \cite [theorem VI.6.1 p.\ 247] {Zygmund_2002}.)} {\it
 Let $a \in \co {} \Zbf$ and consider for a fixed number $q > 1$
 \begin {equation} \label {eq_5.1}
 g(x) \, = \, \sum \sb {m \in \Zbf} a(m) e \sp {i\xx \sb m} \end{equation}
 with $\x \sb {-m} = -\x \sb m$, $\x \sb 1 > 0$, $\x \sb {m + 1} > q\x \sb m$ and $\x \sb m \in \Zbf$ for all $m \in \Zbf \sb +$. Then there is a number $C$ independent of $a$ such that
 \begin {equation*}
 \nm a {} {\lsp \Zbf 1} \, \leq \, C \nm g {} {\Ccal (\Tbf)}. \end{equation*}} \subsection {
 Corollary} \label {cor_5.8} {\it
 Let $\Gcal$ be the set of $g$ appearing in \eqref {eq_5.1} and let $T$ be the linear mapping from $\Gcal$ to $\lsp \Zbf 1$ given by $Tg = a$. Then $T$ has a unique extension, which we also denote by $T$, to the closure of $\Gcal$ in $\Ccal (\Tbf)$ such that
 \begin {equation*}
 \nm {Tg} {} {\lsp \Zbf 1} \, \leq \, C \nm g {} {\Ccal (\Tbf)} \end{equation*}
 where $C$ is independent of $g$.} \section
{Some Lemmata for the Proof of the Main Result} \label {sec_6} \subsection {
 Lemma} \label {lma_6.1} {\it
 Assume that $\Ycal$ is a closed subspace of $\Co \overline n$ such that $\Ycal$ is a subset of $\lkpar \Fcal L \sp 1 \rkpar \lpar \Rn \rpar$. Then there are positive numbers $\alpha$ and $C \sb {\ref {eq_6.01}}$ independent of $\Fcal f \in \Ycal$ such that
\begin {equation} \label {eq_6.01}
\nm {\Fcal f} {} {\Co \overline n} \, \leq \, C \sb {\ref {eq_6.01}} \nm {\Fcal f} {} {\Ccal \lpar \alpha \overline \Bn \, \rpar}. \end{equation}} \subsection*{\it
 Proof:} To simplify notation we write $\overline{\Ccal \sb 0}$, $L \sp 1$ and $c \sb 0$ instead of $\Co \overline n$, $\Lsp n 1$ and $\co {} {\Zbf \sb +}$ respectively. \subsubsection {} \label {sbsbsc_6.1.1}
 Let $\Xcal = \Fcal \sp {-1} \Ycal$. For all $f \in \Xcal$ we have
 \begin {equation} \label {eq_6.02}
 \nm {\Fcal f} {} {\overline{\Ccal \sb 0}} \, \leq \, \nm f {} {L \sp 1} \end{equation}
 and $\Xcal$ is a closed subspace of $L \sp 1$. According to the open mapping theorem there is a number $C \sb {\ref {eq_6.03}}$ independent of $f \in \Xcal$ such that
\begin {equation} \label {eq_6.03}
\nm f {} {L \sp 1} \, \leq \, C \sb {\ref {eq_6.03}}\nm {\Fcal f} {} {\overline{\Ccal \sb 0}}. \end{equation} \subsubsection {}
Let $\eps \sb m$ be a positive number for each $m \in \Zbf \sb +$ such that
\begin{equation*}
\eps \, = \, \sum \sb 1 \sp \infty \eps \sb m \, < \, 1. \end{equation*}
For the purpose of deriving a contradiction we assume that for each choice of positive numbers $\alpha$ and $C$ there is an $\Fcal f \in \Ycal$ such that
\begin{equation*}
C \nm {\Fcal f} {} {\Ccal \left (\alpha \overline \Bn \, \right )} \, < \, \nm {\Fcal f} {} {\overline{\Ccal \sb 0}}. \end{equation*}
As basis for a recursion we choose for $\alpha = \alpha \sb 1 > 0$ and $C = 1/\eps \sb 1$ a function $\Fcal f \sb 1 \in \Xcal$ with $\nm {\Fcal f \sb 1} {} {\overline{\Ccal \sb 0}} = 1$ such that
\begin{equation*}
\nm {\Fcal f \sb 1} {} {\Ccal \left (\alpha \sb 1 \overline \Bn \, \right )} \, < \, \eps \sb 1. \end{equation*}
If $\Fcal f \sb l \in \Ycal$ with $\nm {\Fcal f \sb l} {} {\overline{\Ccal \sb 0}} = 1$ for $l \in \{ 1, \dots , m \}$ as well as $\alpha \sb m > 0$ have been chosen we choose $\alpha \sb {m + 1}$ so that
\begin{equation*}
\sup \lmpar \babs {\lkpar \Fcal f \sb l \rkpar (\x)} : l \in \{ 1, \dots , m \}, \, \axi \geq \alpha \sb {m + 1} \rmpar \, \leq \, \eps \sb {m + 1} \quad \text {and} \quad \alpha \sb {m + 1} > \alpha \sb m. \end{equation*}
By our assumption we can find $\Fcal f \sb {m + 1} \in \Ycal$ with $\nm {\Fcal f \sb {m + 1}} {} {\overline{\Ccal \sb 0}} = 1$ such that
\begin{equation*}
\nm {\Fcal f \sb {m + 1}} {} {\Ccal \left (\alpha \sb {m + 1} \overline \Bn \, \right )} \, < \, \eps \sb {m + 1}. \end{equation*}
We have thus constructed the set $\Phi = \{ f \sb m \in \Xcal : m \in \Zbf \sb + \} \subset L \sp 1$ such that
\begin{align} \label {eq_6.04}
\sup \lmpar \babs {\lkpar \Fcal f \sb l \rkpar (\x)} : l \in \{ 1, \dots , m \}, \, \axi \geq \alpha \sb {m + 1} \rmpar \, & \leq \, \eps \sb {m + 1}, \\
\label {eq_6.05}
\nm {\Fcal f \sb m} {} {\overline{\Ccal \sb 0}} \, &= \, 1 \\ \intertext
{and} \label {eq_6.06}
\nm {\Fcal f \sb {m + 1}} {} {\Ccal \left (\alpha \sb {m + 1} \overline \Bn \, \right )} \, &< \, \eps \sb {m + 1} \end{align}
for each $m \in \Zbf \sb +$ as well as the increasing sequence $\lpar \alpha \sb m \rpar \sb {m \in \Zbf \sb +}$ of positive numbers. \subsubsection {}
For each $m \in \Zbf \sb +$ we choose $b \sb m \in \Rn$ such that $\babs {\lkpar \Fcal f \sb m \rkpar (b \sb m)} = \nm {\Fcal f \sb m} {} {\overline{\Ccal \sb 0}}$ $= 1$. Then $\alpha \sb m < \babs {b \sb m} < \alpha \sb {m + 1}$. Given any $a \in c \sb 0$ such that $a(N + r) = 0$ for some $N \in \Zbf \sb +$ and for all $r \in \Zbf \sb +$ we let $k$ be such that $\babs {a(k) \lkpar \Fcal f \sb k \rkpar \lpar b \sb k \rpar} = \nm a {} {\overline {c \sb 0}}$. Write
\begin{equation*}
\nm a {} {\overline {c \sb 0}} \, = \, \babs {\sum \sb {m = 1} \sp N a(m) \lkpar \Fcal f \sb m \rkpar \lpar b \sb k \rpar \ - \ \sideset {} {\hskip .1mm \sp \prime} \sum \sb {m = 1} \sp N a(m) \lkpar \Fcal f \sb m \rkpar \lpar b \sb k \rpar}, \end{equation*}
where we have omitted the term $a(k) \lkpar \Fcal f \sb k \rkpar \lpar b \sb k \rpar$ from the second sum. We apply the triangle inequality to get
\begin{equation*}
 \nm a {} {\overline{c \sb 0}} \, \leq \, \nm {\sum \sb 1 \sp N a(m) \Fcal f \sb m} {} {\overline{\Ccal \sb 0}} \, + \, \babs {a(k)} \lpar \sum \sb 1 \sp {k - 1} \babs {\lkpar \Fcal f \sb m \rkpar \lpar b \sb k \rpar} + \sum \sb {k + 1} \sp N \babs {\lkpar \Fcal f \sb m \rkpar \lpar b \sb k \rpar} \rpar. \end{equation*}
In the parenthesis we use \eqref {eq_6.04} and \eqref {eq_6.06} on page \pageref {eq_6.06} for the first and second group of terms respectively. (If $k = 1$ or $k = N$ then there is only one group of terms.) We get
\begin{equation*}
\nm a {} {\overline {c \sb 0}} \, < \, \nm {\sum \sb 1 \sp N a(m) \Fcal f \sb m} {} {\overline{\Ccal \sb 0}} + \nm a {} {\overline {c \sb 0}} \, \eps. \end{equation*}
 We have proved that there is a number $C \sb {\ref {eq_6.07}}$ independent of $a \in c \sb 0$ such that
\begin {equation} \label {eq_6.07}
\nm a {} {\overline {c \sb 0}} \, \leq \, C \sb {\ref {eq_6.07}} \nm {\sum \sb {m \in \Zbf \sb +} a(m) \Fcal f \sb m} {} {\overline{\Ccal \sb 0}}. \end{equation}
 \subsubsection {}
There is a vector $b \in \Rn$ such that
\begin{equation*}
\nm {\sum \sb 1 \sp N a(m) \Fcal f \sb m} {} {\overline{\Ccal \sb 0}} \, = \, \babs {\sum \sb 1 \sp N a(m) \lkpar \Fcal f \sb m \rkpar (b)} \, \leq \, \nm a {} {\overline {c \sb 0}} \sum \sb 1 \sp N \babs {\lkpar \Fcal f \sb m\rkpar (b)}. \end{equation*}
Furthermore, there is a unique positive integer $k$ such that $\alpha \sb k \leq |b| < \alpha \sb {k + 1}$. Hence we write
\begin{equation*}
\nm {\sum \sb 1 \sp N a(m) \Fcal f \sb m} {} {\overline{\Ccal \sb 0}} \, \leq \, \nm a {} {\overline {c \sb 0}} \lkpar \lpar \sum \sb 1 \sp {k - 1} + \sum \sb {k + 1} \sp N \rpar \babs {\lkpar \Fcal f \sb m\rkpar (b)} + \babs {\lkpar \Fcal f \sb k \rkpar (b)} \rkpar. \end{equation*}
In the parenthesis we again use \eqref {eq_6.04} and \eqref {eq_6.06} on page \pageref {eq_6.06} for the first and second group of terms respectively. We get
\begin{equation*}
\nm {\sum \sb 1 \sp N a(m) \Fcal f \sb m} {} {\overline{\Ccal \sb 0}} \, < \, \nm a {} {\overline {c \sb 0}} \lpar \eps + 1 \rpar. \end{equation*}
We have proved that there is a number $C \sb {\ref {eq_6.08}}$ independent of $a \in c \sb 0$ such that
\begin {equation} \label {eq_6.08}
\nm {\sum \sb {m \in \Zbf \sb +} a(m) \Fcal f \sb m} {} {\overline{\Ccal \sb 0}} \, \leq \, C \sb {\ref {eq_6.08}} \nm a {} {\overline {c \sb 0}}. \end{equation} \subsubsection {}
For all $a \in c \sb 0$ we have
\begin {equation} \label {eq_6.09}
 \nm a {} {\overline {c \sb 0}} \, \leq \, C \sb {\ref {eq_6.07}} \nm {\sum \sb {m \in \Zbf \sb +} a(m) f \sb m} {} {L \sp 1} \end{equation}
according to \eqref {eq_6.07} and \eqref {eq_6.02} on page \pageref {eq_6.07} and \pageref {eq_6.02} respectively. On the other hand, for all $a \in c \sb 0$ we have
\begin {equation} \label {eq_6.10}
 \nm {\sum \sb {m \in \Zbf \sb +} a(m) f \sb m} {} {L \sp 1} \, \leq \, C \sb {\ref {eq_6.03}} C \sb {\ref {eq_6.08}} \nm a {} {\overline {c \sb 0}} \end{equation}
according to \eqref {eq_6.03} and \eqref {eq_6.08} on page \pageref {eq_6.03} and \pageref {eq_6.08} respectively. \subsubsection {} \label {sbsbsc_6.1.6}
 Let $\Xcal \sb 1$ be the closed linear span of $\Phi$ in $L \sp 1$, and let $\lmpar e \sb m : m \in \Zbf \sb + \rmpar$ be the canonical basis of $\overline {c \sb 0}$. The estimates \eqref {eq_6.09} and \eqref {eq_6.10} give that the mapping $e \sb m \mapsto f \sb m$ can be extended to an isomorphism $\overline {c \sb 0} \longrightarrow \Xcal \sb 1$. But $\Xcal \sb 1$ is a closed subspace of a weakly sequentially complete space and hence $\Xcal \sb 1$ is according to lemma \ref {lma_5.2} on page \pageref {lma_5.2} weakly sequentially complete. Hence $\overline {c \sb 0}$ is weakly sequentially complete. We now invoke \S \ref {sbsbsc_5.1.2} on page \pageref {sbsbsc_5.1.2} to obtain a contradiction. \subsection {}
 In the proof of lemma \ref {lma_6.1} we use only boundedness and linearity \linebreak properties of the Fourier transformation $\Fcal$. Symmetry properties of that linear mapping are not needed for the argument. \subsubsection {\bf
 Proposition} {\it
 Let $\Wcal$ be weakly sequentially complete and let the mapping $T: \Wcal \longrightarrow \Co \overline n$ be bounded, linear and injective. Assume that $\Ycal$ is a closed subspace of $\Co \overline n$ such that $\Ycal$ is a subset of $T\Wcal$. Then there are positive numbers $\alpha$ and $C$ independent of $F \in \Ycal$ such that
\begin{equation*}
\nm F {} {\Co \overline n} \, \leq \, C \nm F {} {\Ccal \left (\alpha \overline \Bn \, \right )}. \end{equation*}} \subsection {
 Lemma} \label {lma_6.3} {\it
 Assume that $\Ycal$ is a closed subspace of $\co \overline \Zn$ such that $\Ycal$ is a subset of $\lkpar \Fcal L \sp 1 \rkpar \lpar \Zn \rpar$. Then there is a positive integer $\alpha$ and a number $C \sb {\ref {eq_6.11}}$ both independent of $\Fcal f \in \Ycal$ such that
\begin {equation} \label {eq_6.11}
\nm {\Fcal f} {} {\co \overline \Zn} \, \leq \, C \sb {\ref {eq_6.11}} \, \sup \lmpar \babs {\lkpar \Fcal f \rkpar (\xi)} : \axi \leq \alpha \rmpar. \end{equation}} \subsection*{\it
 Proof:} The proof is by imitation of the proof of lemma \ref {lma_6.1} on page \pageref {lma_6.1} with some modifications due to the fact that the frequencies are points in $\Zn$. It is given here for the sake of completeness. \par
 To simplify notation we write $L \sp 1$ and $\nm F {} {}$ instead of $L \sp 1 \lpar \Tn \rpar$ and \linebreak $\nm F {} {\co \overline \Zn}$ respectively. Observe that we keep the notation for $\co {} {\Zbf \sb +}$. \subsubsection {}
 Let $\Xcal = \Fcal \sp {-1} \Ycal$. For all $f \in \Xcal$ we have
 \begin {equation} \label {eq_6.12}
 \nm {\Fcal f} {} {} \, \leq \, \nm f {} {L \sp 1} \end{equation}
 and $\Xcal$ is a closed subspace of $L \sp 1$. According to the open mapping theorem there is a number $C \sb {\ref {eq_6.13}}$ independent of $f \in \Xcal$ such that
 \begin {equation} \label {eq_6.13}
 \nm f {} {L \sp 1} \, \leq \, C \sb {\ref {eq_6.13}}\nm {\Fcal f} {} {}. \end{equation} \subsubsection {}
 Let $\eps \sb m$ be a positive number for each $m \in \Zbf \sb +$ such that
 \begin{equation*}
 \eps \, = \, \sum \sb 1 \sp \infty \eps \sb m \, < \, 1. \end{equation*}
 For the purpose of deriving a contradiction we assume that for each choice of $\alpha$ and $C$, where $\alpha$ is a positive integer and $C$ is a positive number, there is an $\Fcal f \in \Ycal$ such that
 \begin{equation*}
 C \, \sup \lmpar \babs {\lkpar \Fcal f \rkpar (\xi)} : \axi \leq \alpha \rmpar \, < \, \nm {\Fcal f} {} {}. \end{equation*}
 As basis for a recursion we choose for a positive integer $\alpha = \alpha \sb 1$ and $C = 1/\eps \sb 1$ a function $\Fcal f \sb 1 \in \Xcal$ with $\nm {\Fcal f \sb 1} {} {} = 1$ such that
 \begin{equation*}
 \sup \lmpar \babs {\lkpar \Fcal f \sb 1 \rkpar (\xi)} : \axi \leq \alpha \sb 1 \rmpar \, < \, \eps \sb 1. \end{equation*}
 If $\Fcal f \sb l \in \Ycal$ with $\nm {\Fcal f \sb l} {} {} = 1$ for $l \in \{ 1, \dots , m \}$ as well as a positive integer $\alpha \sb m$ have been chosen we choose an integer $\alpha \sb {m + 1}$ so that
 \begin{equation*}
 \sup \lmpar \babs {\lkpar \Fcal f \sb l \rkpar (\xi)} : l \in \{ 1, \dots , m \}, \, \axi \geq \alpha \sb {m + 1} \rmpar \leq \eps \sb {m + 1} \quad \text {and} \quad \alpha \sb {m + 1} > \alpha \sb m + 1. \end{equation*}
By our assumption we can find $\Fcal f \sb {m + 1} \in \Ycal$ with $\nm {\Fcal f \sb {m + 1}} {} {} = 1$ such that
\begin{equation*}
 \sup \lmpar \babs {\lkpar \Fcal f \sb {m + 1} \rkpar (\xi)} : \axi \leq \alpha \sb {m + 1} \rmpar \, < \, \eps \sb {m + 1}. \end{equation*}
We have thus constructed the set $\Phi = \{ f \sb m \in \Xcal : m \in \Zbf \sb + \} \subset L \sp 1$ such that
\begin{align} \label {eq_6.14}
\sup \lmpar \babs {\lkpar \Fcal f \sb l \rkpar (\x)} : l \in \{ 1, \dots , m \}, \, \axi \geq \alpha \sb {m + 1} \rmpar \, & \leq \, \eps \sb {m + 1}, \\
\label {eq_6.15}
\nm {\Fcal f \sb m} {} {} \, &= \, 1 \\ \intertext
{and} \label {eq_6.16}
 \sup \lmpar \babs {\lkpar \Fcal f \sb {m + 1} \rkpar (\xi)} : \axi \leq \alpha \sb {m + 1} \rmpar \, &< \, \eps \sb {m + 1} \end{align}
for each $m \in \Zbf \sb +$ as well as the increasing sequence $\lpar \alpha \sb m \rpar \sb {m \in \Zbf \sb +}$ of positive integers. \subsubsection {}
For each $m \in \Zbf \sb +$ we choose $b \sb m \in \Zn$ such that $\babs {\lkpar \Fcal f \sb m \rkpar (b \sb m)} = \nm {\Fcal f \sb m} {} {}$ $= 1$. Then $\alpha \sb m < \babs {b \sb m} < \alpha \sb {m + 1}$. Given any $a \in \co {} {\Zbf \sb +}$ such that $a(N + r) = 0$ for some $N \in \Zbf \sb +$ and for all $r \in \Zbf \sb +$ we let $k$ be such that $\babs {a(k) \lkpar \Fcal f \sb k \rkpar \lpar b \sb k \rpar} = \nm a {} {\co \overline {\Zbf \sb +}}$. Write
 \begin{equation*}
 \nm a {} {\co \overline {\Zbf \sb +}} \, = \, \babs {\sum \sb {m = 1} \sp N a(m) \lkpar \Fcal f \sb m \rkpar \lpar b \sb k \rpar - \sideset {} {\hskip .1mm \sp \prime} \sum \sb {m = 1} \sp N a(m) \lkpar \Fcal f \sb m \rkpar \lpar b \sb k \rpar}, \end{equation*}
 where we have omitted the term $a(k) \lkpar \Fcal f \sb k \rkpar \lpar b \sb k \rpar$ from the second sum. We apply the triangle inequality to get
 \begin{equation*}
 \nm a {} {\co \overline {\Zbf \sb +}} \, \leq \, \nm {\sum \sb 1 \sp N a(m) \Fcal f \sb m} {} {} \, + \, \babs {a(k)} \lkpar \sum \sb 1 \sp {k - 1} \babs {\lkpar \Fcal f \sb m \rkpar \lpar b \sb k \rpar} + \sum \sb {k + 1} \sp N \babs {\lkpar \Fcal f \sb m \rkpar \lpar b \sb k \rpar} \rkpar. \end{equation*}
In the parenthesis we use \eqref {eq_6.14} and \eqref {eq_6.16} on page \pageref {eq_6.14} for the first and second group of terms respectively. (If $k = 1$ or $k = N$ then there is only one group of terms.) We get
 \begin{equation*}
 \nm a {} {\co \overline {\Zbf \sb +}} \, < \, \nm {\sum \sb 1 \sp N a(m) \Fcal f \sb m} {} {} + \nm a {} {\co \overline {\Zbf \sb +}} \eps. \end{equation*}
 We have proved that there is a number $C \sb {\ref {eq_6.17}}$ independent of $a \in \co {} {\Zbf \sb +}$ such that
 \begin {equation} \label {eq_6.17}
\nm a {} {\co \overline {\Zbf \sb +}} \, \leq \, C \sb {\ref {eq_6.17}} \nm {\sum \sb {m \in \Zbf \sb +} a(m) \Fcal f \sb m} {} {}. \end{equation}
 \subsubsection {}
 There is a vector $b \in \Zn$ such that
 \begin{equation*}
 \nm {\sum \sb 1 \sp N a(m) \Fcal f \sb m} {} {} \, = \, \babs {\sum \sb 1 \sp N a(m) \lkpar \Fcal f \sb m \rkpar (b)} \, \leq \, \nm a {} {\co \overline {\Zbf \sb +}} \sum \sb 1 \sp N \babs {\lkpar \Fcal f \sb m\rkpar (b)}. \end{equation*}
 Furthermore, there is a unique positive integer $k$ such that $\alpha \sb k \leq |b| < \alpha \sb {k + 1}$. Hence we write
 \begin{equation*}
 \nm {\sum \sb 1 \sp N a(m) \Fcal f \sb m} {} {} \, \leq \, \nm a {} {\co \overline {\Zbf \sb +}} \lkpar \lpar \sum \sb 1 \sp {k - 1} + \sum \sb {k + 1} \sp N \rpar \babs {\lkpar \Fcal f \sb m \rkpar (b)} + \babs {\lkpar \Fcal f \sb k \rkpar (b)} \rkpar. \end{equation*}
 In the parenthesis we again use \eqref {eq_6.14} and \eqref {eq_6.16} on page \pageref {eq_6.14} for the first and second group of terms respectively. We get
 \begin{equation*}
\nm {\sum \sb 1 \sp N a(m) \Fcal f \sb m} {} {} \, < \, \nm a {} {\co \overline {\Zbf \sb +}} \lpar \eps + 1 \rpar. \end{equation*}
 We have proved that there is a number $C \sb {\ref {eq_6.18}}$ independent of $a \in \co {} {\Zbf \sb +}$ such that
\begin {equation} \label {eq_6.18}
\nm {\sum \sb {m \in \Zbf \sb +} a(m) \Fcal f \sb m} {} {} \, \leq \, C \sb {\ref {eq_6.18}} \nm a {} {\co \overline {\Zbf \sb +}}. \end{equation} \subsubsection {}
 For all $a \in \co {} {\Zbf \sb +}$ we have
 \begin {equation} \label {eq_6.19}
 \nm a {} {\co \overline {\Zbf \sb +}} \, \leq \, C \sb {\ref {eq_6.17}} \nm {\sum \sb {m \in \Zbf \sb +} a(m) f \sb m} {} {L \sp 1} \end{equation}
 according to \eqref {eq_6.17} and \eqref {eq_6.12} on page \pageref {eq_6.17} and \pageref {eq_6.12} respectively. On the other hand, for all $a \in \co {} {\Zbf \sb +}$ we have
 \begin {equation} \label {eq_6.20}
 \nm {\sum \sb {m \in \Zbf \sb +} a(m) f \sb m} {} {L \sp 1} \, \leq \, C \sb {\ref {eq_6.13}} C \sb {\ref {eq_6.18}} \nm a {} {\co \overline {\Zbf \sb +}} \end{equation}
 according to \eqref {eq_6.13} and \eqref {eq_6.18} on page \pageref {eq_6.13} and \pageref {eq_6.18} respectively. \subsubsection {}
 Using \eqref {eq_6.19} and \eqref {eq_6.20} the proof is now concluded in a way similar to the proof of lemma \ref {lma_6.1}. See \S \ref {sbsbsc_6.1.6} on page \pageref {sbsbsc_6.1.6}. \subsection {
 Notation} \label {sbsc_6.4} Let $c \sb k$ be an integer for each $k \in \{ 1, 2, \dots , n \}$ and let $Q = \{ x \in \Rn : c \sb k \leq x \sb k \leq c \sb k + 1, \; k \in \{ 1, 2, \dots , n\}\}$. We form $\Qcal$, the countable collection of all $Q$. In the union representation
\begin{equation*}
 \bigcup \sb {Q \in \Qcal} Q \, = \, \Rbf \sp n \end{equation*}
 the intersection of a pair of terms in the left hand side has Lebesgue measure $0$. If $\beta$ is a positive number we replace $c \sb k$ and $c \sb k + 1$ by $\beta c \sb k$ and $\beta(c \sb k + 1)$ respectively so as to obtain $\beta Q$. For $\beta \neq 1$ the union of $\beta Q$ has the same disjointness property as for $\beta = 1$. \par
 For fixed $\beta > 0$ it is clear that
\begin {equation} \label {eq_6.21}
 \sup \lmpar \babs {x - x \sp \prime} : x, x \sp \prime \in \beta Q, Q \in \Qcal \rmpar \, \leq \, \beta \sqrt n. \end{equation} \subsection {
 Lemma} \label {lma_6.5} {\it
 Let $\Xcal$ be a reflexive subspace of $\Lsp n 1$ of infinite dimension. Then for each choice of positive numbers $\beta$ and $C \sb {\ref {eq_6.22}}$ there is an $f \in \Xcal$ such that
\begin {equation} \label {eq_6.22}
\sum \sb {Q \in \Qcal} \babs {\int \sb {\beta Q} f \ } \, < \, C \sb {\ref {eq_6.22}} \nm f {} {\Lsp n 1}. \end{equation}} \subsection*{\it
 Proof:} For the purpose of deriving a contradiction we assume that there is a choice of positive numbers $\beta$ and $C$ independent of $f \in \Xcal$ such that
\begin{equation*}
 C \nm f {} {\Lsp n 1} \, \leq \, \sum \sb {Q \in \Qcal} \babs {\int \sb {\beta Q} f \ }. \end{equation*}
 Since we also have
\begin{equation*}
\sum \sb {Q \in \Qcal} \babs {\int \sb {\beta Q} f \ } \, \leq \, \sum \sb {Q \in \Qcal} \int \sb {\beta Q} \babs f \, = \, \nm f {} {\Lsp n 1} \end{equation*}
 the mapping $T : \Xcal \longrightarrow \ell \sp 1 (\Zbf \sb +)$ given by
 \begin{equation*}
 [Tf](m) \, = \, \int \sb {\beta Q \sb m} f \ , \quad Q \sb m \in \beta\Qcal \end{equation*}
 is an isomorphism between $\Xcal$ and a closed subspace of $\ell \sp 1 (\Zbf \sb +)$ of infinite dimension.
 According to corollary \ref {cor_5.5} on page \pageref {cor_5.5} this is impossible. This is the contradiction sought for. \subsection {
 Notation} The interval $[0,1[$ is the disjoint union of half-open intervals $[1 - 2 \sp {-k}, 1 - 2 \sp {-k - 1}[$ where $k$ runs through $\Nbf$. The length if each interval is the reciprocal of a dyadic integer. This construction may be transformed to any half-open interval $[a,b[$ using the bijection
\begin{equation*}
\ph : [0,1[ \ \longrightarrow [a,b[ \, , \quad \la \longmapsto (1 - \la)a + \la b. \end{equation*} \par
 Let $N$ be a positive integer. The interval $[-\pi,\pi[$ is the disjoint union of $N$ half-open intervals $[a,b[$ of equal length $\beta = 2\pi/N$. For each such interval \linebreak $[a,b[$ we apply the construction using subintervals of $[0,1[$ of length of a \linebreak reciprocal dyadic integer and the bijection $\ph$. As $[a,b[$ runs through finitely many subintervals of $\Tbf = [-\pi,\pi[$ we obtain a countable disjoint union \linebreak representing that interval. By definition, a term in this union representation is a {\it $\beta$-admissable interval}. \par
Let $\beta = 2\pi/N$ be given for some positive integer $N$. For each factor $\Tbf$ in the cartesian product $\Tn$ we pick a $\beta$-admissable interval and form the cartesian product $R$ of those $n$ intervals. We also form $\Rcal \sb \beta$, the countable collection of all $R$. In the union representation
\begin{equation*}
 \bigcup \sb {R \in \Rcal \sb \beta} R \, = \, \Tbf \sp n \end{equation*}
 the intersection of a pair of terms in the left hand side has Lebesgue measure $0$. \par
 For fixed $\beta > 0$ it is clear that there is a positive number $C$ such that
\begin {equation} \label {eq_6.23}
 \sup \lmpar \babs {x - x \sp \prime} : x, x \sp \prime \in R, R \in \Rcal \sb \beta \rmpar \, \leq \, \beta C. \end{equation} \subsection {
 Lemma} \label {lma_6.7} {\it
 Let $\Xcal$ be a reflexive subspace of $L \sp 1 (\Tn)$ of infinite dimension. Then for each choice of positive numbers $\beta$ and $C \sb {\ref {eq_6.24}}$, where $\beta = 2\pi/N$ for $N \in \Zbf \sb +$, there is an $f \in \Xcal$ such that
 \begin {equation} \label {eq_6.24}
\sum \sb {R \in \Rcal \sb \beta} \babs {\int \sb R f \ } \, < \, C \sb {\ref {eq_6.24}} \nm f {} {L \sp 1 \lpar \Tn \rpar}. \end{equation}} \subsection*{\it
 Proof:} We imitate the proof of lemma \ref {lma_6.5} on page \pageref {lma_6.5} whereby
 \begin{equation*}
 [Tf](m) \, = \, \int \sb {R \sb m} f \ , \quad R \sb m \in \Rcal \sb \beta. \end{equation*} \section
{Proof of Theorem \ref {thm_4.1}} \label {sec_7} \subsection {\it
 Proof of theorem {\rm \ref {thm_4.1}} on page {\rm \pageref {thm_4.1}} in the case $X = \Rn$:} \label {sbsc_7.1} \subsubsection {}
 For each $f \in L \sp 1 = \Lsp n 1$, for each $m \in \Zbf \sb +$ and for each $\beta > 0$ we have with the notation from \S \ref {sbsc_6.4} on page \pageref {sbsc_6.4}
 \begin{multline*}
 \babs {\int \sb {\beta Q \sb m} e \sp {-i\xx} f(x) \dx} \, \leq \, \int \sb {\beta Q \sb m} \babs {e \sp {-i\xx} - e \sp {-iq \sb m \x}} \babs {f(x)} \dx \, + \\
 + \, \babs {e \sp {-iq \sb m \x} \int \sb {\beta Q \sb m} f \ } \, \leq \, \int \sb {\beta Q \sb m} \babs {x - q \sb m} \axi \babs {f(x)} \dx + \babs {\int \sb {\beta Q \sb m} f \ } \end{multline*}
 for any $q \sb m \in Q \sb m$. Summing with respect to $m$, taking $\sup$ with respect to $\x$ and invoking \eqref {eq_6.21} gives
\begin {equation} \label {eq_7.1}
 \nm {\Fcal f} {} {\Ccal \left (\alpha \overline \Bn \, \right )} \leq \sup \sb {\axi \leq \alpha} \sum \sb 1 \sp \infty \babs {\int \sb {\beta Q \sb m} e \sp {-i\xx} f(x) \dx} \leq \alpha \beta \sqrt n \nm f {} {L \sp 1} + \sum \sb 1 \sp \infty \babs {\int \sb {\beta Q \sb m} f }. \end{equation} \subsubsection {}
 Assume that $\Ycal$ fulfills the assumptions of the theorem. As in the proof of lemma \ref {lma_6.1} (cf.\ \S \ref {sbsbsc_6.1.1} on page \pageref {sbsbsc_6.1.1}) the space $\Xcal = \Fcal \sp {-1} \Ycal$ is according to the open mapping theorem isomorphic to $\Ycal$. We have
 \begin {equation} \label {eq_7.2}
 \nm f {} {L \sp 1} \, \leq \, C \sb {\ref {eq_6.03}} \nm {\Fcal f} {} {\Co \overline n} \, \leq \, C \sb {\ref {eq_6.03}} C \sb {\ref {eq_6.01}} \nm {\Fcal f} {} {\Ccal \left (\alpha \overline \Bn \, \right )} \end{equation}
 where we have used lemma \ref {lma_6.1} on page \pageref {lma_6.1} in the second inequality. Collecting the estimates \eqref {eq_7.1} and \eqref {eq_7.2} gives that there is a number $C \sb {\ref {eq_7.3}}$ independent of $f$ such that
\begin {equation} \label {eq_7.3}
 \nm f {} {L \sp 1} \, \leq \, C \sb {\ref {eq_7.3}} \lkpar \alpha \beta \sqrt n \nm f {} {L \sp 1} + \sum \sb 1 \sp \infty \babs {\int \sb {\beta Q \sb m} f \ } \rkpar. \end{equation} \subsubsection {}
 By assumption, $\Ycal$ is reflexive. Hence $\Xcal$ is reflexive. For the purpose of deriving a contradiction we now assume that $\Xcal$ is of infinite dimension. Then the assumptions of lemma \ref {lma_6.5} on page \pageref {lma_6.5} are fulfilled, and hence for 
\begin{equation*}
 \beta \, < \, \frac 1 {2 \, C \sb {\ref {eq_7.3}} \, \alpha \sqrt n} \quad \text {and} \quad C \sb {\ref {eq_6.22}} \, = \, \frac 1 {2C \sb {\ref {eq_7.3}}} \end{equation*}
 there is an $f \in \Xcal$ such that
\begin{equation*}
 \nm f {} {L \sp 1} \leq C \sb {\ref {eq_7.3}} \lkpar \alpha \beta \sqrt n \nm f {} {L \sp 1} + \sum \sb 1 \sp \infty \babs {\int \sb {\beta Q \sb m} f \ } \rkpar < \frac 1 2 \nm f {} {L \sp 1} + \frac 1 2 \nm f {} {L \sp 1} = \nm f {} {L \sp 1}. \end{equation*}
 This is the contradiction sought for. \subsection {\it
 Remarks on the proof of theorem {\rm \ref {thm_4.1}} on page {\rm \pageref {thm_4.1}} in the case $X = \Tn$:} \label {sbsc_7.2} The proof is by imitation of the proof in the case $X = \Rn$ which was just completed. Lemmata \ref {lma_6.1} and \ref {lma_6.5} on page \pageref {lma_6.1} and \pageref {lma_6.5} respectively are replaced by lemmata \ref {lma_6.3} and \ref {lma_6.7} on page \pageref {lma_6.3} and \pageref {lma_6.7} respectively. Inequality \eqref {eq_6.21} on page \pageref {eq_6.21} is replaced by inequality \eqref {eq_6.23} on page \pageref {eq_6.23}. 
\section
{Examples of Closed Non-Reflexive Subspaces} \label {sec_8} \subsection {}
 The examples provided here are based on Karlander's idea in \cite [p.\ 312] {Karlander_1997} using lacunary trigonometric series. \subsection {
 The Case $X = \Rn$.} \label {sbsc_8.2} Let $H$ be any {\it positive} function in $\Lsp n 1 \cap \Co \overline n$. We say that the function $F$ belongs to the space $\Ycal$ if and only if there is an $a \in \lsp {\Zbf \sb +} 1$ such that
 \begin {equation*}
 F(\x) \, = \, \lkpar Ta \rkpar (\x) \, = \, H(\x) \sum \sb {k \in \Zbf \sb+} a(k) e \sp {i2 \sp k \x \sb 1}, \quad \x \, = \, (\x \sb 1, \dots , \x \sb n) \in \Rn. \end{equation*}
 Then $F$ is the Fourier transform of an $\Lsp n 1$-function. \subsubsection {}
 We have
 \begin {equation*}
 \nm {Ta} {} {\Co \overline n} \, \leq \, \nm H {} {\Co \overline n} \nm a {} {\lsp {\Zbf \sb +} 1} \end{equation*}
 and so $T$ is a bounded linear bijection from $\lsp {\Zbf \sb +} 1$ to $\Ycal$. Assume that we can show that $\Ycal$ is closed. Then, according to the open mapping theorem $T \sp {-1}$ is a bounded linear bijection from $\Ycal$ to $\lsp {\Zbf \sb +} 1$, and we may conclude that $\Ycal$ and $\lsp {\Zbf \sb +} 1$ are isomorhic. In particular, $\Ycal$ is not reflexive. This shows that we may obtain closed subspaces of $\lkpar \Fcal L \sp 1 \rkpar \lpar \Rn \rpar$ of infinite dimension if we drop the reflexivity requirement. \subsubsection {}
 We now show that $\Ycal$ is closed. \par
 Assume that $F \sb m$ is in $\Ycal$ for each $m \in \Zbf \sb +$ and that $F \sb m$ converges to $F$ in $\Co \overline n$ as $m \to \infty$. If $G \sb m = F \sb m/H$ and if $K$ is a compact set then there is a function $G \in \Ccal \lpar \Rn \rpar$ such that $G \sb m$ converges to $G$ in $\Ccal (K)$ as $m \to \infty$. But
 \begin{equation*}
 G \sb m (\x) \, = \, \sum \sb {k \in \Zbf \sb+} a \sb m (k) e \sp {i2 \sp k \x \sb 1} \end{equation*}
 for some $a \sb m \in \lsp {\Zbf \sb +} 1$. We now invoke corollary \ref {cor_5.8} on page \pageref {cor_5.8} to conclude that there is a function $a \in \lsp {\Zbf \sb +} 1$ such that $a \sb m$ converges to $a$ in $\lsp {\Zbf \sb +} 1$ as $m \to \infty$. For fixed $\x \in \Rn$ we have
 \begin{equation*}
 F(\x) \, = \, \lim \sb {m \to \infty} H(\xi) \sum \sb {k \in \Zbf \sb+} a \sb m (k) e \sp {i2 \sp k \x \sb 1} \, = \, H(\xi) \sum \sb {k \in \Zbf \sb+} a(k) e \sp {i2 \sp k \x \sb 1}. \end{equation*}
 We have proved that $F \in \Ycal$ and hence $\Ycal$ is closed. \subsection {
 The Case $X = \Tn$.} Let $H$ be any {\it positive} function in $\lsp \Zn 1$. We say that the function $F$ belongs to the space $\Ycal$ if and only if there is an $a \in \lsp {\Zbf \sb +} 1$ such that
 \begin {equation*}
 F(\x) \, = \, H(\x) \sum \sb {k \in \Zbf \sb+} a(k) e \sp {i2 \sp k \x \sb 1}, \quad \x \, = \, (\x \sb 1, \dots , \x \sb n) \in \Zn. \end{equation*}
 Then $F$ is the Fourier transform of an $L \sp 1 \lpar \Tn \rpar$-function. \par
 In a way similar to the argument in \S \ref {sbsc_8.2} on page \pageref {sbsc_8.2} one may prove that $\Ycal$ is a closed subspace of $\co \overline \Zn$ such that $\Ycal$ is a subset of $\lkpar \Fcal L \sp 1 \rkpar \lpar \Zn \rpar$ and such that $\Ycal$ is isomorphic to $\lsp {\Zbf \sb +} 1$. Thus, also in this case we may obtain closed subspaces $\Ycal$ of $\co \overline \Zn$ of infinite dimension such that $\Ycal$ is a subset of $\lkpar \Fcal L \sp 1 \rkpar \lpar \Zn \rpar$ if we drop the reflexivity requirement.
  \end{document}